\documentclass{amsart}

\usepackage{xypic}
 \input xy
\xyoption{all}
\usepackage{epsfig}
\usepackage{amsthm}
\usepackage{amssymb}
\usepackage{amsmath}
\usepackage{amscd}
\usepackage{color}
\usepackage{mathabx}

\newcommand{\bg}{\begin{equation}}
\newcommand{\ed}{\end{equation}}
\newcommand{\bga}{\begin{eqnarray}}
\newcommand{\eda}{\end{eqnarray}}
\newcommand{\pf}{\textbf{Proof:\ }}

\def\cbdu{\par{\raggedleft$\Box$\par}}

\newtheorem {Theorem}  {Theorem}

\numberwithin{Theorem}{section}

\newtheorem {Lemma}[Theorem]  {Lemma}

\theoremstyle{definition}
\newtheorem{Definition}[Theorem]{Definition}
\theoremstyle{remark}
\newtheorem{Remark}[Theorem]{\bf Remark}

\expandafter\chardef\csname pre amssym.def
at\endcsname=\the\catcode`\@ \catcode`\@=11
\def\undefine#1{\let#1\undefined}
\def\newsymbol#1#2#3#4#5{\let\next@\relax
 \ifnum#2=\@ne\let\next@\msafam@\else
 \ifnum#2=\tw@\let\next@\msbfam@\fi\fi
 \mathchardef#1="#3\next@#4#5}
\def\mathhexbox@#1#2#3{\relax
 \ifmmode\mathpalette{}{\m@th\mathchar"#1#2#3}%
 \else\leavevmode\hbox{$\m@th\mathchar"#1#2#3$}\fi}
\def\hexnumber@#1{\ifcase#1 0\or 1\or 2\or 3\or 4\or 5\or 6\or 7\or 8\or
 9\or A\or B\or C\or D\or E\or F\fi}

\font\teneufm=eufm10 \font\seveneufm=eufm7 \font\fiveeufm=eufm5
\newfam\eufmfam
\textfont\eufmfam=\teneufm \scriptfont\eufmfam=\seveneufm
\scriptscriptfont\eufmfam=\fiveeufm

\catcode`\@=\csname pre amssym.def at\endcsname

\newcounter{remark}
\usepackage{color}

\newcommand{\e}{\epsilon}

\renewcommand{\th}{\theta}

\newcommand{\R}{\mathbf{R}}
\renewcommand{\div}{\mbox{div}}

\def  \R   {{\mathbb R}}

\def  \12  {{\frac{1}{2}}}

\def\build#1_#2^#3{\mathrel{\mathop{\kern 0pt#1}\limits_{#2}^{#3}}}

\begin{document}

\title[Stability of Quasi-Geostrophic Equation]{Existence and stability of steady-state solutions to the Quasi-Geostrophic equations in $\mathbb R^2$}

\author[Mimi Dai]{ Mimi Dai}
\address{Department of Mathematics, University of Illinois at Chicago, Chicago, IL 60607,USA}
\email{mdai@uic.edu}





\begin{abstract}
 We consider the stationary Quasi-Geostrophic equation in the whole space $\mathbb R^2$ driven by a force $f$. Under  certain smallness assumptions of $f$, we establish the existence of solutions with finite $L^2$ norm. This solution is unique among all solutions with finite energy. The unique solution $\Theta$ is also shown to be stable in the sense: any solution of the evolutionary Quasi-Geostrophic equation driven by $f$ and starting with finite energy, will return to $\Theta$.
\end{abstract}

\maketitle

\section{Introduction}
In this paper we consider the existence and stability of stationary solutions to the two dimensional active scalar equation

\begin{equation}\begin{split}\label{QG}
\theta_t+u\cdot\nabla \theta+\kappa\Lambda^\alpha\theta =f,\\
u=R^\perp\theta,
\end{split}
\end{equation}
in $\R^2\times(0, T)$ with $T>0$, where $\alpha\in[1,2)$, $\kappa>0$, $\Lambda=\sqrt{-\Delta}$ is the Zygmund operator, and
\bg\notag
R^\perp\theta=\Lambda^{-1}(-\partial_2\theta,\partial_1\theta).
\ed
 The scalar function $\theta$ represents the potential temperature and the vector function $u$ represents the fluid velocity. The force $f$ is independent of time. If $\theta_t\equiv 0$, (\ref{QG}) reduces to a stationary equation:
\begin{equation}\begin{split}\label{SQG}
U\cdot\nabla \Theta+\kappa\Lambda^\alpha\Theta =f,\\
U=R^\perp\Theta.
\end{split}
\end{equation}

Equation (\ref{QG}) with $\alpha=1$ is a fundamental model of the surface quasi-geostrophic equation (SQG) \cite{CMT, Pe}. It describes the evolution of the surface temperature field in a rapidly rotating and stably stratified fluid with potential velocity. As pointed out in \cite{CMT}, this equation attracts interest of  scientists and mathematicians due to two major reasons: it is a fundamental model for the actual geophysical flows with applications in atmosphere and oceanography study; from the mathematical point of view, the behavior of strongly nonlinear solutions to (\ref{QG}) with $\kappa=0, f=0$ in 2D and the behavior of potentially singular solutions to the Euler's equation in 3D are strikingly analogous which has been justified both analytically and numerically. For literature the readers are refereed to \cite{CCW, CMT, CW,  CC, DKSV, DKV, Don, Pe} and the references therein. 

When $\alpha=1$, equation (\ref{QG}) is usually referred as critical SQG, 
although it is an open problem whether a dramatic change in the behavior of solutions occurs for the case of dissipation power less than $1$. The global regularity problem of the critical SQG equation has been very challenging due to the balance of the nonlinear term and the dissipative term
in (\ref{QG}). This problem for the unforced critical SQG is resolved now by Kieslev, Nazarov and Volberg \cite{KNV}, Caffarelli and Vasseur \cite{CaV}, Kieslev and Nazarov \cite{KN09, KN10} and Constantin and Vicol \cite{CV} independently, using different sophisticated methods.  

The long time behaviors of the solutions to the critical SQG equation have been studied in \cite{CD, CTVan, CTV, DD, FPV, NS, SS03, SS05}. Estimates of the decay rates are obtained for mild solutions, regular solutions and weak solutions in \cite{SS05, DD}, \cite{NS}, and \cite{SS03}, respectively. 
In the forced case in \cite{CTVan}, with a special class of time independent force, the long time average behavior of viscosity solutions has been addressed and the absence of anomalous dissipation was obtained. 
 In \cite{FPV}, the authors proved that the linear instability implies the nonlinear instability in the energy norm in the periodic setting (on domain $\mathbb T^2$). Recently, Constantin, Tarfulea and Vicol \cite{CTV} and Cheskidov and Dai \cite{CD} 
proved the existence of a compact global attractor in $H^1(\mathbb T^2)$ and $L^2(\mathbb T^2)$, respectively, with different assumptions on $f$.
In particular, in \cite{CD}, the assumption is that $f$ is solely in $L^p(\mathbb T^2)$ for some $p>2$. A crucial estimate is that any initial data $\theta_0\in L^2(\mathbb T^2)$ yields weak (viscosity) solutions with bounded norm in $L^\infty$,  by using the De Giorgi iteration method. We point out that the boundedness in
$L^\infty$ also plays an important rule in the stability analysis of the current paper.


In this paper, we study the existence and uniqueness of solutions to (\ref{SQG}) with finite energy ($L^2$ norm) in the whole space $\mathbb R^2$, and the nonlinear stability of such solution. On a domain where the Poincar\'e's type of inequality holds it is expected that the steady state has finite $L^2$ norm; and when the force is small, any evolutionary solution with finite $L^2$ norm converges to the unique steady state. But the Poincar\'e's type of inequality is not available in the whole space $\mathbb R^2$, where the situation is more delicate. Thus techniques other than the standard energy method are applied to study these problems. To establish the existence of solutions to (\ref{SQG}) with finite energy, inspired by the work of Bjorland and Schonbek \cite{BS} for Navier-Stokes equation, we make use of the fact:
if $\Phi(t,x)$ is a fundamental solution for the operator $\partial_t+\Lambda^\alpha=0$, then $\int_0^\infty\Phi(t,\cdot)dt$ is a fundamental solution for the operator $\Lambda^\alpha=0$. Therefore, the $L^2$ estimate of the steady state may be obtained through the fast decay estimates of solutions to the corresponding evolutionary equation. In the mean time, the Fourier splitting method (see \cite{Sch1, Sch2}) is at hand to establish the decay estimates of the evolutionary solutions in the whole space $\mathbb R^2$. Once the finite energy estimate holds, the smallness assumption on the force yields the uniqueness of the steady state. To study the stability of the steady state, we combine the method of generalized energy estimates for low and high frequency parts of the solutions and the Fourier splitting method.

In contrast to the NSE as studied in \cite{BS}, the SQG equation \eqref{QG} has a weaker dissipative term $\Lambda^\alpha$ with $\alpha\in[1,2)$. 
This presents a serious obstacle to establish the stability of the steady state.
However, thanks to the fact that a viscosity solution is bounded in $L^\infty$ (see \cite{CDmodes, CD}), this difficulty can be overcome. Thus, the stability is obtained in the sense: viscosity solutions of \eqref{QG} converge to the steady state $\Theta$.

Before stating the main results, we first recall the definitions of weak solution and viscosity solution.

\begin{Definition}
A weak solution to \eqref{QG} is a function $\th \in C_{\mathrm{w}}([0,T];L^2(\mathbb{R}^2))$ satisfying that , for any $\phi\in C_0^\infty(\mathbb R^2\times (0,T))$,
\begin{equation}\notag
-\int_0^T(\theta, \phi_t)\, dt-\int_0^T(u\theta, \nabla\phi)\, dt+\kappa\int_0^T(\Lambda^{\frac{1}{2}}\theta, \Lambda^{\frac{1}{2}}\phi)\, dt
=(\theta_0, \phi(x,0))+\int_0^T(f, \phi)\, dt.
\end{equation}
\end{Definition}

\begin{Definition}
A weak solution $\th(t)$ of (\ref{QG}) on $[0,T]$ is called a viscosity solution if there exist sequences $\e_n \to 0$ and $\th_n(t)$ satisfying
\begin{equation}\begin{split}\label{VQG}
\frac{\partial\theta_n}{\partial t}+u_n\cdot\nabla \theta_n+\kappa\Lambda^\alpha\theta_n + \e_n \Delta \th_n=f,\\
u_n=R^\perp\theta_n,
\end{split}
\end{equation}
such that $\th_n \to \theta$ in $C_\mathrm{w}([0,T];L^2)$. 
\end{Definition}

Standard arguments imply that for any initial data $\th_0 \in L^2$ there exists a viscosity solution $\th(t)$ of \eqref{QG} with $\alpha=1$ on $[0,\infty)$ with $\th(0)=\th_0$ (see \cite{CC}, for example). The same result holds for the subcritical SQG equation with $1<\alpha<2$. It was proved in \cite{CDmodes, CD} that a viscosity solution to (\ref{QG}) is bounded in $L^\infty$ provided $f\in L^p$ for some $p>\frac2\alpha$. With this auxiliary estimate in $L^\infty$, we are able to prove the nonlinear stability of the steady state $\Theta$. 

The main results for the critical SQG ($\alpha=1$) and the subcritical SQG with $\alpha\in(1,2)$ are stated separately as follows.

\begin{Theorem}\label{thm:ex}
Suppose $f\in X=W^{\frac{1}{2}, 4}\cap L^\infty$ satisfy the assumption:\\
(A) $\hat f(\xi)=0$ for almost every $|\xi|<\rho_0$, for some constant $\rho_0>0$.\\
Then there exist a constant $M>0$ and a constant $C(\rho_0,\kappa, M)$ so that if $\|f\|_X\leq C(\rho_0,\kappa, M)$, 
system (\ref{SQG}) with $\alpha=1$ has a weak solution $\Theta\in H^{\frac12}$, in the sense that for any $\phi\in\mathcal V$,
\bg\notag
-(U\Theta,\nabla\phi)+\kappa(\Lambda^{\frac1{2}} U,\Lambda^{\frac1{2}}\Theta)
=(f,\phi).
\ed
which satisfies 
\begin{equation}\label{energy-Theta}
\|\Theta\|_2\leq M, \quad \|\Lambda^{\frac1{2}}\Theta\|_2\leq\kappa^{-1}\|f\|_X.
\end{equation}
Moreover, this solution is unique among all solutions satisfying \eqref{energy-Theta}.
\end{Theorem}

\begin{Theorem}\label{thm:stability}
There exists a constant $C(\kappa)$ such that if $\|f\|_X\leq C(\kappa)$,  the steady state $\Theta$ obtained in Theorem \ref{thm:ex} is nonlinearly stable in the sense: let $w_0\in L^2$ and $\theta$ be a viscosity solution of (\ref{QG}) with $\alpha=1$ and initial data $\theta_0=w_0+\Theta$ which satisfies, for any $T>0$,
$$
\theta\in L^\infty(0,T;L^2)\cap L^2(0,T; \dot H^{\frac{1}{2}})
$$
then
\[\lim_{t\to\infty}\|\theta(t)-\Theta\|_2=0.\]
\end{Theorem}

\begin{Remark}
The constant $C(\kappa)$ in Theorem \ref{thm:stability} may be even smaller than the constant $C(\rho_0,\kappa, M)$ in Theorem \ref{thm:ex}. The smallness of $\|f\|_X\leq C(\kappa)$ will imply the smallness of $\Theta$ in $H^1$ by Lemma  \ref{le:higher}. It is not known whether the smallness assumption on the force in Theorem \ref{thm:ex} or Theorem \ref{thm:stability} can be removed or not. To apply the techniques presented in this paper, the smallness assumption is crucial. It is also worth to mention that after the work of \cite{BS}, the authors and their collaborators proved the stability of the steady state of 3D NSE with rougher external force in \cite{BBIS}. Namely, the assumption (A) was relaxed, but the smallness of $f$ was still required.
\end{Remark}

In the subcritical case $\alpha\in(1,2)$, we establish that:

\begin{Theorem}\label{thm:sub}
Let $\alpha\in(1,2)$. Suppose $f\in X=W^{1-\frac{\alpha}{2}, 4}\cap L^{\frac{4}{\alpha-1}}$ satisfying the assumption A.
Then there exist a constant $M>0$ and another constant $C(\rho_0,\kappa, M)$ so that if $\|f\|_X\leq C(\rho_0,\kappa, M)$ the following hold:\\
(i) System (\ref{SQG}) has a weak solution $\Theta$, in the sense that for any $\phi\in\mathcal V$,
\bg\notag
-(U\Theta,\nabla\phi)+\kappa(\Lambda^{\frac{\alpha}{2}} \Theta,\Lambda^{\frac{\alpha}{2}}\phi)
=(f,\phi).
\ed
And $\Theta\in H^{\frac{\alpha}{2}}$, with 
\[\|\Theta\|_2\leq M, \quad \|\Lambda^{\frac{\alpha}{2}}\Theta\|_2\leq\kappa^{-1}\|f\|_X.\]
Moreover, this solution is unique among all solutions with finite norm in $L^2$.\\
(ii) The steady state $\Theta$ is nonlinearly stable in the sense: let $w_0\in L^2$ and $\theta$ be a viscosity solution of the SQG equation (\ref{QG}) with initial data $\theta_0=w_0+\Theta$ which satisfies, for any $T>0$,
$$
\theta\in L^\infty(0,T;L^2)\cap L^2(0,T; \dot H^{\frac{\alpha}{2}})
$$
then
\[\lim_{t\to\infty}\|\theta(t)-\Theta\|_2=0.\]
\end{Theorem}

In these theorems and through the paper, we adopt the notations, $\|\cdot\|_p=\|\cdot\|_{L^p}$, $(f, g)=\int_{\mathbb R^2}fgdx$ and $\mathcal V=\{\phi\in C^\infty_0|\nabla\cdot \phi=0\}$.

We shall prove the results for the critical SQG ($\alpha=1$), that is, Theorem \ref{thm:ex} and Theorem \ref{thm:stability}. For the subcritical SQG with $\alpha\in(1,2)$, a similar analysis scheme with slight modification will yield the analogous result. Thus the proof of Theorem \ref{thm:sub} will be omitted.

The paper is organized as follows: in section 2 we introduce the proof scheme and  give some preliminary estimates; section 3 is devoted to the proof of the existence of steady state with finite energy, that is, Theorem \ref{thm:ex}; and section 4 is devoted to the proof of the stability of the steady state, that is,
Theorem \ref{thm:stability}.

\medskip

\section{Preliminaries} 
\label{sec:pre}

\subsection{Outline of the analysis scheme}
The nontrivial part in the study of the stationary equation (\ref{SQG}) is to establish $\|\Theta\|_2<M$ due to the fact that a Poincar\'e type of inequality is not available in the whole space $\R^2$. Inspired by the work of \cite{BS} for the Navier-Stokes equation, we plan to achieve the goal in the following. Consider
\begin{equation}\begin{split}\label{QGf}
\frac{\partial\tilde\theta}{\partial t}+U\cdot\nabla \tilde\theta+\kappa\Lambda^\alpha\tilde\theta =0\\
\tilde\theta(0)=f, \ \ \ U=R^\perp\Theta.
\end{split}
\end{equation}
Formally, if $\tilde\theta$ solves (\ref{QGf}), $\tilde\Theta=\int_0^\infty\tilde\theta(t)dt$ solves
\begin{equation}\begin{split}\label{QGst}
U\cdot\nabla \tilde\Theta+\kappa\Lambda^\alpha\tilde\Theta =f\\
U=R^\perp\Theta.
\end{split}
\end{equation}
Since this equation is linear for a fixed $U$, 
solutions are unique and hence $\tilde\Theta=\Theta$. Due to the integral Minkowski's inequality, 
\bg\notag
\|\Theta\|_2=\|\int_0^\infty\tilde\theta(t)dt\|_2\leq\int_0^\infty\|\tilde\theta(t)\|_2dt,
\ed
one can see $\Theta\in L^2$ if $\|\tilde\theta(t)\|_2\leq C(1+t)^{-\gamma}$ with $\gamma>1$. Thus it is crucial to establish fast decay for $\tilde \theta$ which will be addressed in Section \ref{sec:pf1}.

\subsection{Estimates on the operator $\Lambda$}  We recall some standard estimates which will be used often through the paper.

\begin{Lemma}\label{le:phi}
Let $\Phi(t)=e^{-\kappa t\Lambda^\alpha f}$.
If $f$ satisfies Assumption (A), then
\begin{equation}\notag
\|\Lambda^\nu\Phi(t)\|_p\leq e^{-\kappa \rho_0^\alpha t}\|\Lambda^\nu f\|_p
\end{equation}
for any $1\leq p\leq\infty$ and $\nu\geq 0$, where we adopt the convention $\Lambda^0v\equiv v$.
\end{Lemma}
For a proof of Lemma \ref{le:phi}, see \cite{BS}.


\begin{Lemma}\label{le:Sob} \cite{Le}(Sobolev type inequality)
Let $2<p<\infty$ and $\nu=1-\frac{2}{p}$. There exists a constant $C\geq 0$ such that if $v\in\mathcal S'$ and $\hat v$ is a function, then 
\bg\notag
\|v\|_p\leq C\|\Lambda^\nu v\|_2.
\ed
\end{Lemma}
\begin{Lemma}\label{le:tem-vel}
Let $1<p<\infty$. There exists a constant $C_p$ depending only on $p$ such that
\bg\notag
\|\Lambda^\nu u(t)\|_p\leq C_p\|\Lambda^\nu \theta(t)\|_p
\ed
for all $\nu\geq 0$, $t\geq 0$.
\end{Lemma}



\bigskip

\section{Existence of Steady State with Finite $L^2$ Norm}
\label{sec:pf1}

As outlined in Section \ref{sec:pre}, the fast decay of solution $\tilde\theta$ to (\ref{QGf}) as $\|\tilde\theta(t)\|_2\leq C(1+t)^{-\gamma}$ with $\gamma>1$ plays an important rule. Due to a result in \cite{CW}, the best decay one can expect for the critical case is that
(since $U$ does not decay, it is worse than the real SQG in decay rate)
\bg\notag
\|\tilde\theta(t)\|_2\leq C(1+t)^{-1}
\ed
through a standard Fourier splitting method (c.f. \cite{Sch1, Sch2}), where the initial data prevents a faster decay. 
To obtain a faster decay, we measure the difference 
\bg\notag
\beta=\tilde\theta-\Phi
\ed
with 
$ \Phi=e^{-\kappa t\Lambda^\alpha}f$, and $\beta$ satisfies 
\begin{equation}\begin{split}\label{QGb}
\beta_t+U\cdot\nabla \beta+\kappa\Lambda^\alpha\beta =-U\cdot\nabla\Phi\\
\beta(0)=0, \ \ \ U=R^\perp\Theta.
\end{split}
\end{equation}
We expect that $\|\beta(t)\|_2$ decays as $(1+t)^{-\frac{2}{\alpha}}$ (see the proof for $\alpha=1$ in Section \ref{sec:pf1}). By assumption (A), the decay of $\Phi$ is fast enough such that $\tilde\theta$ will decay as fast as $\beta$.

We first study the solutions for the two sequences of approximating equations:
\begin{equation}\label{eq:app1}
\begin{split}
U^i\cdot\nabla \Theta^{i+1}+\kappa\Lambda\Theta^{i+1} =f\\
U^i=R^\perp\Theta^i,
\end{split}
\end{equation}
\begin{equation}\begin{split}\label{eq:appt1}
\beta^{i+1}_t+U^i\cdot\nabla \beta^{i+1}+\kappa\Lambda\beta^{i+1} =-U^i\cdot\nabla\Phi\\
\beta^{i+1}(0)=0, \ \ \ U^i=R^\perp\Theta^i.
\end{split}
\end{equation}
For a fixed function $U^i\in H^{\frac{1}{2}}$, we solve the two systems recursively to find approximating solutions for (\ref{SQG}) and (\ref{QGb}) with $\alpha=1$. With uniform estimates on the approximating solutions, the convergence procedure will yield the estimate of $\|\Theta\|_2<M$.

\medskip

\subsection{Existence of Solutions to the Approximating Systems}

\begin{Lemma}\label{thm:appex1}
Let $U^i\in H^{\frac{1}{2}}$ and $f\in X$. There exists a unique weak solution $\Theta^{i+1}$
to (\ref{eq:app1}) in the sense that for any $\phi\in\mathcal V$,
\bg\label{eq:weak1}
(U^i\Theta^{i+1},\nabla\phi)+\kappa(\Lambda^{\frac{1}{2}}\Theta^{i+1},\Lambda^{\frac{1}{2}}\phi)=(f,\phi).
\ed
Moreover, this solution satisfies
\bg\label{thetai1}
\|\Lambda^{\frac{1}{2}}\Theta^{i+1}\|_2\leq C\kappa^{-1}\|f\|_X.
\ed
\end{Lemma}
For the linear equation (\ref{eq:app1}), a standard Galerkin method gives the existence of solutions $\Theta^{i+1}$ which
 satisfy
\begin{equation}\notag
\kappa\int_{\R^2}|\Lambda^{\frac{1}{2}}\Theta^{i+1}|^2dx\leq\int_{\R^2}f\Theta^{i+1}dx
\leq C\|\Lambda^{-\frac{1}{2}}f\|_2\|\Lambda^{\frac{1}{2}}\Theta^{i+1}\|_2
\end{equation}
with an absolute constant $C$. 


\begin{Lemma}\label{thm:apptex1}
Let $U^i\in H^{\frac{1}{2}}$ with
\bg\label{ineq:ui1}
\|\Lambda^{\frac{1}{2}}U^{i}\|_2\leq C\kappa^{-1}\|f\|_X
\ed
and $f$ satisfy assumption A.
 There exists a unique weak solution $\beta^{i+1}\in L^\infty(\R^+, L^2)\cap L^2(\R^+, \dot H^{\frac{1}{2}})$
to (\ref{eq:appt1}) in the sense that for any $\phi\in C^1(\R^+;\mathcal V)$,
\bg\notag
\begin{split}
-(\beta^{i+1},\phi_t)-(U^i\beta^{i+1},\nabla\phi)+\kappa(\Lambda^{\frac1{2}}\beta^{i+1},\Lambda^{\frac1{2}}\phi)\\
=-(U^i\cdot\nabla\Phi,\phi)\\
\beta^{i+1}(x,0)=0, \ \ \ U^i=R^\perp\Theta^i.
\end{split}
\ed
Moreover, it satisfies
\bg\label{est:betai1}
\sup_{t>0}\left(\|\beta^{i+1}(t)\|_2^2+\kappa\int_0^t\|\Lambda^{\frac{1}{2}}\beta^{i+1}(s)\|_2^2ds\right)\leq C\rho_0^{-1}\kappa^{-4}\|f\|_X^4.
\ed
\end{Lemma}
\pf
Applying Galerkin method yields a sequence of approximating smooth solutions $\beta^{i+1}_m$ which satisfy
\begin{align}\notag
&\frac{1}{2}\frac{d}{dt}\int_{\R^2}|\beta^{i+1}_m|^2dx+\kappa\int_{\R^2}|\Lambda^{\frac{1}{2}}\beta^{i+1}_m|^2dx\\
&=-\int_{\R^2}(U^i\cdot\nabla\Phi)\beta^{i+1}_mdx\notag\\
&\leq \frac{C}{\kappa}\|\Lambda^{\frac{1}{2}}(U^i\Phi)\|_2^2+\frac{\kappa}{2}\|\Lambda^{\frac{1}{2}}\beta^{i+1}_m\|_2^2\notag.
\end{align}
It follows
\bg\label{ineq:bit}
\frac{d}{dt}\int_{\R^2}|\beta^{i+1}_m|^2dx+\kappa\int_{\R^2}|\Lambda^{\frac{1}{2}}\beta^{i+1}_m|^2dx\leq\frac{C}{\kappa}\|\Lambda^{\frac{1}{2}}(U^i\Phi)\|_2^2.
\ed
Applying the H\"older inequality and Lemma \ref{le:phi}, we have
\begin{align}\notag
\|\Lambda^{\frac{1}{2}}(U^i\Phi)\|_2^2&\leq C\|\Lambda^{\frac{1}{2}}U^i\Phi\|_2^2+
C\|U^i\Lambda^{\frac{1}{2}}\Phi\|_2^2\\
&\leq C\|\Lambda^{\frac{1}{2}}U^i\|_2^2\|\Phi\|_\infty^2+
C\|U^i\|_4^2\|\Lambda^{\frac{1}{2}}\Phi\|_4^2\notag\\
&\leq C\|\Lambda^{\frac{1}{2}}U^i\|_2^2\|\Phi\|_\infty^2+
C\|\Lambda^{\frac{1}{2}}U^i\|_2^2\|\Lambda^{\frac{1}{2}}\Phi\|_4^2\notag\\
&\leq C\kappa^{-2}e^{-2\kappa\rho_0t}\left(\|f\|_\infty^2+\|\Lambda^{\frac{1}{2}}f\|_4^2\right)\|f\|_X^2\notag.
\end{align}
Integrating over the time interval $[0,t]$ for (\ref{ineq:bit}) yields
\begin{align}\notag
\int_{\R^2}|\beta^{i+1}_m(t)|^2dx+\kappa\int_0^t\int_{\R^2}|\Lambda^{\frac{1}{2}}\beta^{i+1}_m|^2dx&\leq C\kappa^{-3}\|f\|_X^4\int_0^te^{-2\kappa\rho_0s}ds\\
&\leq C\rho_0^{-1}\kappa^{-4}\|f\|_X^4\notag.
\end{align}
The inequality (\ref{est:betai1}) follows by taking the limit $m\to\infty$.

\cbdu

\subsection{Decay of $\beta^i$}
\label{sec:betadec}
In this subsection we establish the decay for $\beta^i$ by using the Fourier Splitting method.
 The rigorous proof can be carried through by working on a sequence of the Galerkin approximating solutions to (\ref{eq:appt1}). To be abbreviate, we work on the regular solutions in a formal way. 

\begin{Lemma}\label{le:Fbeta}
Let $\beta^{i+1}$ be the solution of (\ref{eq:appt1}) given by Lemma \ref{thm:apptex1}.
The following estimate holds
\bg\label{est:betadecay}
\|\beta^{i+1}(t)\|_2^2
\leq C(\rho_0,\kappa)\|U^i\|_2^4\|f\|_X^4(1+t)^{-4}
\ed
with the constant $C(\rho_0,\kappa)$ depending on $\rho_0$ and $\kappa$.
\end{Lemma}
\pf 
The first step is to establish 
\bg\label{est:Fbeta}
|\hat\beta^{i+1}(t)|\leq C|\xi|\|U^i\|_2\left(\int_0^t\|\beta^{i+1}\|_2(s)ds+\rho_0^{-1}\kappa^{-1}\|f\|_2\right).
\ed
Indeed, taking Fourier transform of (\ref{eq:appt1}) yields
\bg\notag
\hat\beta^{i+1}_t+\kappa|i\xi|\hat\beta^{i+1}=-\left(\mathcal F(U^i\cdot\nabla\beta^{i+1})+\mathcal F(U^i\cdot\nabla\Phi)\right).
\ed
Since $\beta^{i+1}(0)=0$, we have
\bg\notag
\hat\beta^{i+1}=-\int_0^te^{-\kappa|\xi|(t-s)}\left(\mathcal F(U^i\cdot\nabla\beta^{i+1})+\mathcal F(U^i\cdot\nabla\Phi)\right)(s)ds.
\ed
Due to the fact that $\nabla\cdot U^i=0$, it follows by the Young's inequality and the Plancherel theorem
\begin{equation}\notag
\begin{split}
|\hat\beta^{i+1}|&\leq C|\xi|\int_0^t|\mathcal F(U^i\beta^{i+1})(s)|+|\mathcal F(U^i\Phi)(s)|ds\\
&\leq C|\xi|\int_0^t\|(U^i\beta^{i+1})(s)\|_{L^1}+\|(U^i\Phi)(s)\|_{L^1}ds\\
&\leq C|\xi|\|U^i\|_2\int_0^t\|\beta^{i+1}(s)\|_2+\|\Phi(s)\|_2ds\\
&\leq C|\xi|\|U^i\|_2\int_0^t\|\beta^{i+1}(s)\|_2+e^{-\kappa\rho_0s}\|f\|_2ds\\
&\leq C|\xi|\|U^i\|_2\left(\int_0^t\|\beta^{i+1}\|_2(s)ds+\rho_0^{-1}\kappa^{-1}\|f\|_2\right).
\end{split}
\end{equation} 

Then we claim that, for any $l>5$, $\beta^{i+1}$ satisfies
\begin{equation}\label{est:betaik}
\begin{split}
&\frac{d}{dt}\left[(1+t)^l\|\beta^{i+1}\|_2^2\right]\\
&\leq C(l,\rho_0,\kappa)(1+t)^{l-5}\|U^i\|_2^2\left(\int_0^t\|\beta^{i+1}(s)\|_2ds+\|f\|_2\right)^2\\
&+C\kappa^{-3}\|f\|_X^2\left(\|\Phi\|_\infty^2+\|\Lambda^{\frac{1}{2}}\Phi\|_4^2\right)(1+t)^{l},
\end{split}
\end{equation}
with the constant $C(l,\rho_0,\kappa)=C\kappa(\frac{l}{\kappa})^5(1+\rho_0^{-1}\kappa^{-1})^2$ for an absolute constant $C$. 

In order to prove it, multiplying (\ref{eq:appt1}) by $\beta^{i+1}$ and integrating over $\R^2$ yields
\begin{align}\label{ineq:betai2}
\frac{d}{dt}\|\beta^{i+1}\|_2^2+\kappa\|\Lambda^{\frac{1}{2}}\beta^{i+1}\|_2^2
&\leq C\kappa^{-1}\left(\|\Lambda^{\frac{1}{2}}U^i\|_2^2\|\Phi\|_\infty^2+\|\Lambda^{\frac{1}{2}}U^i\|_2^2\|\Lambda^{\frac{1}{2}}\Phi\|_4^2\right)\\
&\leq C\kappa^{-3}\|f\|_X^2\left(\|\Phi\|_\infty^2+\|\Lambda^{\frac{1}{2}}\Phi\|_4^2\right)\notag.
\end{align}
Denote $S(t)$ as a sphere in $\R^2$ with center at the origin and radius $R(t)=\frac{l}{\kappa}(1+t)^{-1}$.  We decompose the viscous term in frequency space into two parts, corresponding to $S(t)$ and its complementary $S(t)^c$, respectively. 
Applying the Plancherel theorem yields
\begin{align}\notag
-\kappa\|\Lambda^{\frac{1}{2}}\beta^{i+1}\|_2^2&\leq-\kappa\int_{S(t)^c}|\xi||\hat\beta^{i+1}|^2d\xi\\
&\leq-\kappa R\int_{S(t)^c}|\hat\beta^{i+1}|^2d\xi\notag\\
&\leq-\kappa R\|\hat\beta^{i+1}\|_2^2+\kappa R\int_{S(t)}|\hat\beta^{i+1}|^2d\xi\notag.
\end{align}
Combining (\ref{ineq:betai2}) and the last inequality it follows that
\begin{align}\notag
\frac{d}{dt}\|\beta^{i+1}\|_2^2+\kappa R\|\beta^{i+1}\|_2^2
&\leq \kappa R\int_{S(t)}|\hat\beta^{i+1}|^2d\xi+C\kappa^{-3}\|f\|_X^2\left(\|\Phi\|_\infty^2+\|\Lambda^{\frac{1}{2}}\Phi\|_4^2\right)\notag.
\end{align}
By  (\ref{est:Fbeta}) we have
\begin{align}\notag
\int_{S(t)}|\hat\beta^{i+1}|^2d\xi&\leq C\|U^i\|_2^2\left(\int_0^t\|\beta^{i+1}(s)\|_2ds+\rho_0^{-1}\kappa^{-1}\|f\|_2\right)^2\int_{S(t)}|\xi|^2d\xi\\
&\leq C\|U^i\|_2^2(1+\rho_0^{-1}\kappa^{-1})^2\left(\int_0^t\|\beta^{i+1}(s)\|_2ds+\|f\|_2\right)^2R^4\notag.
\end{align}
Thus,
\begin{align}\notag
\frac{d}{dt}\|\beta^{i+1}\|_2^2+\kappa R\|\beta^{i+1}\|_2^2
&\leq C\kappa\|U^i\|_2^2(1+\rho_0^{-1}\kappa^{-1})^2\left(\int_0^t\|\beta^{i+1}(s)\|_2ds+\|f\|_2\right)^2R^5\notag\\
&+C\kappa^{-3}\|f\|_X^2\left(\|\Phi\|_\infty^2+\|\Lambda^{\frac{1}{2}}\Phi\|_4^2\right)\notag
\end{align}
Multiplying the last inequality by the time factor $(1+t)^l$, it follows 
\begin{align}\notag
\frac{d}{dt}\left[(1+t)^l\|\beta^{i+1}\|_2^2\right]
&\leq C(l,\rho_0,\kappa)\|U^i\|_2^2(1+t)^{l-5}\left(\int_0^t\|\beta^{i+1}(s)\|_2ds+\|f\|_2\right)^2\notag\\
&+C\kappa^{-3}\|f\|_X^2\left(\|\Phi\|_\infty^2+\|\Lambda^{\frac{1}{2}}\Phi\|_4^2\right)(1+t)^l\notag
\end{align}
with $C(l,\rho_0,\kappa)=C\kappa(\frac{l}{\kappa})^5(1+\rho_0^{-1}\kappa^{-1})^2$ for an absolute constant $C$, which concludes the argument of the claim.


Combining (\ref{est:betai1}) and (\ref{est:betaik}) and using the fact $\|f\|_2\leq \|f\|_X$ yields
\begin{align}\notag
&\frac{d}{dt}\left[(1+t)^l\|\beta^{i+1}\|_2^2\right]\\
&\leq C(l,\rho_0,\kappa)(1+t)^{l-5}\|U^i\|_2^2\left(t\rho_0^{-\frac{1}{2}}\kappa^{-2}\|f\|_X^2+\|f\|_2\right)^2\notag\\
&+C\kappa^{-3}\|f\|_X^2\left(\|\Phi\|_\infty^2+\|\Lambda^{\frac{1}{2}}\Phi\|_4^2\right)(1+t)^{l}\notag\\
&\leq C(l,\rho_0,\kappa)\|U^i\|_2^2\left(1+\|f\|_X^2\right)\|f\|_X^2(1+t)^{l-3}\notag\\
&+C\kappa^{-3}\|f\|_X^2\left(\|\Phi\|_\infty^2+\|\Lambda^{\frac{1}{2}}\Phi\|_4^2\right)(1+t)^{l}\notag.
\end{align}
Integrating in time and applying Lemma \ref{le:phi}, it follows that
\begin{align}\notag
&(1+t)^l\|\beta^{i+1}(t)\|_2^2\\
&\leq C(l,\rho_0,\kappa)\|U^i\|_2^2\left(1+\|f\|_X^2\right)\|f\|_X^2\int_0^t(1+s)^{l-3}ds\notag\\
&+C\kappa^{-3}\|f\|_X^2\int_0^t\left(\|\Phi(s)\|_\infty^2+\|\Lambda^{\frac{1}{2}}\Phi(s)\|_4^2\right)(1+s)^{l}ds\notag\\
&\leq C(l,\rho_0,\kappa)\|U^i\|_2^2\left(1+\|f\|_X^2\right)\|f\|_X^2\left[(1+t)^{l-2}-1\right]\notag\\
&+C\kappa^{-3}\|f\|_X^2\left(\|f\|_\infty^2+\|\Lambda^{\frac{1}{2}}f\|_4^2\right)\int_0^te^{-2\kappa\rho_0s}(1+s)^{l}ds\notag.
\end{align}
One can verify that 
\[\int_0^te^{-2\kappa\rho_0s}(1+s)^{l}\, ds\sim (1+t)^le^{-2\kappa\rho_0 t}\]
which is bounded for all $t>0$ since $l, \rho_0, \kappa>0$. It infers
\[
C\kappa^{-3}\|f\|_X^2\left(\|f\|_\infty^2+\|\Lambda^{\frac{1}{2}}f\|_4^2\right)\int_0^te^{-2\kappa\rho_0s}(1+s)^{l}ds\leq C(l,\rho_0,\kappa) \|f\|_X^4.
\] 
Thus, we obtain the preliminary decay estimate
\bg\label{est:betadec0}
\|\beta^{i+1}(t)\|_2^2\leq C(l,\rho_0,\kappa)\|U^i\|_2^2\left(1+\|f\|_X^2\right)\|f\|_X^2(1+t)^{-2}.
\ed
Combining (\ref{est:betaik}) and (\ref{est:betadec0}), we repeat the procedure as above and obtain that
\bg\notag
\|\beta^{i+1}(t)\|_2^2\leq C(l,\rho_0,\kappa)\|U^i\|_2^4\|f\|_X^4(1+t)^{-4}.
\ed
\cbdu

\subsection{Proof of $\Theta^{i+1}=\int_0^\infty\tilde\theta^{i+1}dt$ solving (\ref{eq:app1})}
As pointed out in Section \ref{sec:pre}, formally $\Theta^{i+1}=\int_0^\infty\tilde\theta^{i+1}dt$ solves the approximating system (\ref{eq:app1}) with $\tilde\theta^{i+1}=\beta^{i+1}+\Phi$. In this subsection, we prove this statement rigorously and show that $\Theta^{i+1}$ is uniformly bounded in $L^2$.
\begin{Lemma}\label{le:theta-l2}
Let $\beta^{i+1}$ be the solution of (\ref{eq:appt1}) given by Theorem \ref{thm:apptex1}. Then it satisfies $\int_0^\infty\tilde\theta^{i+1}(t)dt=\int_0^\infty(\beta^{i+1}+\Phi)dt\in L^2$. As a consequence, $\int_0^\infty\tilde\theta^{i+1}(t)dt$ is finite $a.e.$ in $\R^2$.
\end{Lemma}
\pf
For each $i$ we define the sequence $\left\{\Theta^{i+1}_n\right\}$ as
\bg\notag
\Theta^{i+1}_n=\int_0^n\tilde\theta^{i+1}(t)dt.
\ed
We have $\Theta^{i+1}_n\in L^2$ since $\tilde\theta^{i+1}(t)\in L^2$ for almost every $t$ by Lemma \ref{thm:apptex1}. Applying Minkowski's inequality, Lemma \ref{le:phi} and (\ref{est:betadecay}) yields
\begin{align}\notag
\|\Theta^{i+1}_n\|_2&\leq \int_0^n\|\tilde\theta^{i+1}(t)\|_2dt\\
&\leq \int_0^n\|\beta^{i+1}(t)\|_2dt+\int_0^n\|\Phi(t)\|_2dt\notag\\
&\leq C(\rho_0,\kappa)\|U^i\|_2^2\|f\|_X^2+\|f\|_X\notag.
\end{align}
Similarly, the Minkowski's inequality implies 
\bg\notag
\|\Theta^{i+1}_{n+1}-\Theta^{i+1}_{n}\|_2\leq \int_n^{n+1}\|\tilde\theta^{i+1}(t)\|_2dt.
\ed
It is known that $\int_0^{\infty}\|\tilde\theta^{i+1}(t)\|_2dt$ is finite by Lemmas \ref{le:phi} and  \ref{le:Fbeta}. Thus, 
\[\int_n^{n+1}\|\tilde\theta^{i+1}(t)\|_2dt\to 0 \qquad \mbox { as } n\to \infty. \]
We infer that the sequence $\left\{\Theta^{i+1}_n\right\}$ is Cauchy in $L^2$ and hence it converges to a limit $\int_0^\infty\tilde\theta^{i+1}(t)dt\in L^2$. 
\cbdu

\begin{Lemma}\label{le:Theta-sol}
We have $\Theta^{i+1}=\int_0^\infty\tilde\theta^{i+1}(t)dt$, that is, $\int_0^\infty\tilde\theta^{i+1}(t)dt$ solves (\ref{eq:app1}) with $U^i$ satisfying (\ref{ineq:ui1}) and $f$ satisfying the assumption A.
\end{Lemma}
\pf
Multiplying (\ref{eq:appt1}) by $\phi\in\mathcal V$, integrating over space and time, we infer that from $\beta^{i+1}=\tilde\theta^{i+1}-\Phi$
\bg\notag
\int_0^n\left(\frac{d}{dt}(\tilde\theta^{i+1}(t),\phi)+(\Lambda^{\frac12}(U^i\tilde\theta^{i+1}),\Lambda^{\frac12}\phi)\right)dt=-\kappa\int_0^n(\Lambda^{\frac{1}{2}}\tilde\theta^{i+1},\Lambda^{\frac{1}{2}}\phi)dt.
\ed
Noticing that $\tilde\theta^{i+1}(0)=\beta^{i+1}(0)-\Phi(0)=f$, exchanging the order of the integration yields
\bg\label{app-n}
(\tilde\theta^{i+1}(n),\phi)+(\Lambda^{\frac12}(U^i\Theta^{i+1}_n),\Lambda^{\frac12}\phi)=-\kappa(\Lambda^{\frac12}\Theta^{i+1}_n,\Lambda^{\frac{1}{2}}\phi)+(f,\phi).
\ed
Since $\tilde\theta^{i+1}(n)=\beta^{i+1}(n)+\Phi(n)$, it follows that the first term on the left hand side tends to $0$ as $n\to\infty$ by (\ref{est:betadecay}).
On the other hand, denote $\tilde\Theta^{i+1}=\int_0^\infty\tilde\theta^{i+1}(t)dt$, and we have
\begin{align}\notag
|(\Lambda^{\frac12}[U^i(\Theta^{i+1}_n-\tilde\Theta^{i+1})],\Lambda^{\frac12}\phi)|&\leq C\|U^i\|_4\|\Theta^{i+1}_n-\tilde\Theta^{i+1}\|_2\|\nabla\phi\|_4\\
&\leq C\|\Lambda^{\frac12} U^i\|_2\|\Theta^{i+1}_n-\tilde\Theta^{i+1}\|_2\|\nabla\phi\|_4\notag
\end{align}
\bg\notag
|(\Lambda^{\frac{1}{2}}(\Theta^{i+1}_n-\tilde\Theta^{i+1}),\Lambda^{\frac12}\phi)|\leq
C\|\Theta^{i+1}_n-\tilde\Theta^{i+1}\|_2\|\Lambda\phi\|_2.
\ed
Lemma \ref{le:theta-l2} implies that the sequence $\left\{\Theta^{i+1}_n\right\}$ converges strongly to $\tilde\Theta^{i+1}=\int_0^\infty\tilde\theta^{i+1}(t)dt$ in $L^2$. Thus, we conclude from the last two inequalities that 
\begin{align}\notag
&|(U^i\cdot\nabla(\Theta^{i+1}_n-\tilde\Theta^{i+1}),\phi)|\to 0, \ \ \ n\to\infty\\
&|(\Lambda^{\frac{1}{2}}(\Theta^{i+1}_n-\tilde\Theta^{i+1}),\Lambda^{\frac{1}{2}}\phi)|\to 0, \ \ \ n\to\infty\notag.
\end{align}
Therefore taking the limit $n\to\infty$ in (\ref{app-n}) yields that $\tilde\Theta^{i+1}$ is a weak solution of (\ref{eq:app1}). The uniqueness in Lemma \ref{thm:appex1} implies $\tilde \Theta^{i+1}=\Theta^{i+1}$, which completes the proof of the lemma.
\cbdu

\subsection{Finishing the proof of Theorem \ref{thm:ex}} 
\label{sec:conv}
In this subsection we show that the approximating solutions $\Theta^i$ to (\ref{eq:app1}) converge to a limit function $\Theta$ which solves the stationary equation (\ref{SQG}) with certain smallness assumptions on $f$. Additionally, we establish a uniqueness result in the sense stated in Theorem \ref{thm:ex}.

\begin{Lemma}\label{le:l2bd}
Let $\Theta^{i+1}$ be the solution of (\ref{eq:app1}) given by Theorem \ref{thm:appex1} with $U^i$ and $f$ satisfying the assumptions of the theorem. Let $M>0$ be a constant such that
$\|U^i\|_2\leq M$. Then there exists a constant $C(\rho_0,\kappa, M)$ such that if $\|f\|_X\leq C(\rho_0,\kappa, M)$ then $\|\Theta^{i+1}\|_2\leq M$.
\end{Lemma}

\pf
By Lemma \ref{le:Theta-sol},  Lemma \ref{le:phi} and (\ref{est:betadecay}) we infer that
\begin{align}\label{const-determine}
\|\Theta^{i+1}\|_2&\leq \int_0^\infty\|\theta^{i+1}(t)\|_2dt\\
&\leq \int_0^\infty\|\beta^{i+1}(t)\|_2dt+ \int_0^\infty\|\Phi(t)\|_2dt\notag\\
&\leq C(\rho_0,\kappa)\|U^i\|_2^2\|f\|_X^2+\|f\|_X\notag
\end{align}
Take $M$ satisfying the quadratic equation
\bg\notag
C(\rho_0,\kappa)M^2Z^2+Z=M
\ed
with a positive root $Z=\frac{-1+\sqrt{1+4C(\rho_0,\kappa)M^3}}{2C(\rho_0,\kappa)M^2}=: C(\rho_0,\kappa, M)$.
It follows from (\ref{const-determine}) that if $\|f\|_X\leq C(\rho_0,\kappa, M)$, then $\|\Theta^{i+1}\|_2\leq M$.
By the boundedness of Riesz transform on $L^2$, it follows that $\|U^{i+1}\|_2\leq M$.
\cbdu

The following lemma establishes the estimates for higher order derivatives of the sequence $\Theta^{i+1}$.

\begin{Lemma}\label{le:higher}
Let $\Theta^{i+1}$ be the solution of (\ref{eq:app1}) given by Lemma \ref{thm:appex1} with $U^i$ and $f$ satisfying the assumptions of the lemma. Assume additionally that $\Theta^i\in H^{\frac{3}{2}}$ with
\bg\label{ass:1}
\|\Lambda\Theta^i\|_2\leq C\kappa^{-1}\|f\|_X, \ \ \ \|\Lambda^{\frac{3}{2}}\Theta^i\|_2\leq C\kappa^{-1}\|f\|_X.
\ed
 Then there exists a constant $C(\kappa, M)$  such that if $\|f\|_X\leq C(\kappa, M)$ then $\Theta^{i+1}\in H^{\frac{3}{2}}$ satisfying
\bg\notag
\|\Lambda\Theta^{i+1}\|_2\leq C\kappa^{-1}\|f\|_X, \ \ \ \|\Lambda^{\frac{3}{2}}\Theta^{i+1}\|_2\leq C\kappa^{-1}\|f\|_X.
\ed
The constant $C(\kappa, M)$ depends on $\kappa$ and $M$. 
\end{Lemma}
\pf Formally, multiplying (\ref{eq:app1}) by $\Lambda\Theta^{i+1}$ (to make it rigorous we can use the Galerkin approximating solutions as test function here) yields
\begin{align}\notag
\kappa\|\Lambda\Theta^{i+1}\|_2^2&=-(U^i\cdot\nabla\Theta^{i+1},\Lambda\Theta^{i+1})+
(f,\Lambda\Theta^{i+1})\\
&\leq \|U^i\|_\infty\|\Lambda\Theta^{i+1}\|_2^2+C\kappa^{-1}\|f\|_2^2+\frac{\kappa}{4}\|\Lambda\Theta^{i+1}\|_2^2\notag.
\end{align}
Followed from the classical Gagliardo-Nirenberg inequalities \cite{EG} by complex interpolation (see also \cite{KP}), and by the boundedness of Riesz transform we have 
\begin{align}\label{ineq:ui-infty}
\|U^i\|_\infty&\leq C\|\Lambda^{\frac{3}{2}}U^i\|_2^{\frac{2}{3}}\|U^i\|_2^{\frac{1}{3}}\\
&\leq C\|\Lambda^{\frac{3}{2}}\Theta^i\|_2^{\frac{2}{3}}\|\Theta^i\|_2^{\frac{1}{3}}\notag.
\end{align}
Combining the last two inequalities, the assumption (\ref{ass:1}) and Lemma \ref{le:l2bd} gives that
\begin{align}\notag
\|\Lambda\Theta^{i+1}\|_2^2
&\leq C\kappa^{-1}\|\Lambda^{\frac{3}{2}}\Theta^i\|_2^{\frac{2}{3}}\|\Theta^i\|_2^{\frac{1}{3}}\|\Lambda\Theta^{i+1}\|_2^2+C\kappa^{-2}\|f\|_2^2\notag\\
&\leq C\kappa^{-\frac{5}{3}}M^{\frac{1}{3}}\|f\|_X^{\frac{2}{3}}\|\Lambda\Theta^{i+1}\|_2^2+C\kappa^{-2}\|f\|_2^2\notag.
\end{align}
We choose $C(\kappa, M)$ such that $\|f\|_X\leq C(\kappa, M)$ and $C\kappa^{-\frac{5}{3}}M^{\frac{1}{3}}\|f\|_X^{\frac{2}{3}}\leq 1/2$. Thus,
\bg\notag
\|\Lambda\Theta^{i+1}\|_2^2\leq C\kappa^{-2}\|f\|_2^2.
\ed
Similarly, we multiply (\ref{eq:app1}) formally by $\Lambda^2\Theta^{i+1}$ and infer that
\begin{align}\notag
\kappa\|\Lambda^{\frac{3}{2}}\Theta^{i+1}\|_2^2&=-(U^i\cdot\nabla\Theta^{i+1},\Lambda^2\Theta^{i+1})+
(f,\Lambda^2\Theta^{i+1})\\
&\leq |(\Lambda^{\frac{1}{2}}U^i\cdot\nabla\Theta^{i+1},\Lambda^{\frac{3}{2}}\Theta^{i+1})|+
|(U^i\cdot\Lambda^{\frac{1}{2}}\nabla\Theta^{i+1},\Lambda^{\frac{3}{2}}\Theta^{i+1})|\notag\\
&+C\kappa^{-1}\|\Lambda^{\frac{1}{2}}f\|_2^2+\frac{\kappa}{4}\|\Lambda^{\frac{3}{2}}\Theta^{i+1}\|_2^2\notag\\
&\leq C\kappa^{-1}\|\Lambda^{\frac{1}{2}}U^i\cdot\nabla\Theta^{i+1}\|_2^2+
\|U^i\|_\infty\|\Lambda^{\frac{3}{2}}\Theta^{i+1}\|_2^2\notag\\
&+C\kappa^{-1}\|\Lambda^{\frac{1}{2}}f\|_2^2+\frac{\kappa}{2}\|\Lambda^{\frac{3}{2}}\Theta^{i+1}\|_2^2\notag.
\end{align}
By Lemma \ref{le:Sob} and the boundedness of Riesz Transform it follows
\begin{align}\notag
\|\Lambda^{\frac{1}{2}}U^i\cdot\nabla\Theta^{i+1}\|_2^2
&\leq C\|\Lambda^{\frac{1}{2}}U^i\|_4^2\|\nabla\Theta^{i+1}\|_4^2\\
&\leq C\|\Lambda\Theta^i\|_2^2\|\Lambda^{\frac{3}{2}}\Theta^{i+1}\|_2^2\notag.
\end{align}
Combining (\ref{ineq:ui-infty}) and the last two inequalities yields
\begin{align}\notag
\|\Lambda^{\frac{3}{2}}\Theta^{i+1}\|_2^2&\leq C\kappa^{-2}\|\Lambda^{\frac{1}{2}}f\|_2^2+C\kappa^{-2}\|\Lambda\Theta^{i}\|_2^2\|\Lambda^{\frac{3}{2}}\Theta^{i+1}\|_2^2\notag\\
&+C\kappa^{-1}\|\Lambda^{\frac{3}{2}}\Theta^{i}\|_2^{\frac{2}{3}}\|\Theta^{i}\|_2^{\frac{1}{3}}\|\Lambda^{\frac{3}{2}}\Theta^{i+1}\|_2^2\notag\\
&\leq C\kappa^{-2}\|\Lambda^{\frac{1}{2}}f\|_2^2
+C\left(\kappa^{-4}\|f\|_X^2+\kappa^{-\frac{5}{3}}M^{\frac{1}{3}}\|f\|_X^{\frac{2}{3}}\right)\|\Lambda^{\frac{3}{2}}\Theta^{i+1}\|_2^2\notag.
\end{align}
Thus if we choose $C(\kappa, M)$ such that $\|f\|_X\leq C(\kappa, M)$ and 
$$C\left(\kappa^{-4}\|f\|_X^2+\kappa^{-\frac{5}{3}}M^{\frac{1}{3}}\|f\|_X^{\frac{2}{3}}\right)\leq 1/2$$ 
then
\bg\notag
\|\Lambda^{\frac{3}{2}}\Theta^{i+1}\|_2^2\leq C\kappa^{-2}\|f\|_X^2.
\ed

\cbdu

We are now ready to prove the first main theorem.\\
\textbf{Proof of Theorem \ref{thm:ex}:}
Choose $\Theta^0\in H^{3/2}$ such that $\|\Theta^0\|_2\leq M$ and $\|\Lambda^{\gamma}\Theta^0\|_2\leq \kappa^{-1}\|f\|_X$ with $\gamma=\frac{1}{2},1,\frac{3}{2}$. By the boundedness of Riesz transform on Sobolev spaces, we have $U^0=R^\perp\Theta^0\in H^{3/2}$, $\|U^0\|_2\leq M$ and $\|\Lambda^{\gamma}U^0\|_2\leq \kappa^{-1}\|f\|_X$ with $\gamma=\frac{1}{2},1,\frac{3}{2}$.

Starting with $U^0$, we solve (\ref{eq:app1}) recursively by Lemma \ref{thm:appex1} and obtain a sequence $\left\{\Theta^i\right\}_{i=0}^\infty$ with $U^i=R^\perp\Theta^i$. It satisfies that $\|\Lambda^{\gamma}\Theta^i\|_2\leq \kappa^{-1}\|f\|_X$ with $\gamma=\frac{1}{2},1,\frac{3}{2}$ by Lemma \ref{thm:appex1} and Lemma \ref{le:higher}, and $\|\Theta^i\|_2\leq M$ uniformly by Lemma \ref{le:l2bd}. In the following we show that the sequence $\left\{\Theta^i\right\}_{i=0}^\infty$ is Cauchy in $\dot H^{1/2}$ and it converges to a limit function $\Theta$.

Let $Y^{i+1}=\Theta^{i+1}-\Theta^i$. Thus $Y^{i+1}$ solves the equation
\bg\notag
U^i\cdot\nabla Y^{i+1}+R^\perp Y^i\cdot\nabla\Theta^i+\kappa\Lambda Y^{i+1}=0.
\ed
Since $Y^{i+1}\in H^{1/2}$, we can multiply the above equation by $Y^{i+1}$, integrate by parts and we have, by noticing that $\nabla\cdot U^i=\nabla\cdot R^\perp Y^i=0$
\begin{align}\notag
\kappa\|\Lambda^{\frac{1}{2}}Y^{i+1}\|_2^2&=-(\nabla\cdot(R^\perp Y^i\Theta^i), Y^{i+1})\\
&\leq C\kappa^{-1}\|\Lambda^{\frac{1}{2}}(R^\perp Y^{i}\Theta^i)\|_2^2+\frac{\kappa}{2}\|\Lambda^{\frac{1}{2}}Y^{i+1}\|_2^2\notag.
\end{align}
Thus by the boundedness of Riesz transform, Gagliardo-Nirenberg inequalities and Lemmas \ref{le:l2bd} and \ref{le:higher} we infer that
\begin{align}\label{eq:recur}
\|\Lambda^{\frac{1}{2}}Y^{i+1}\|_2^2
&\leq C\kappa^{-2}\|\Lambda^{\frac{1}{2}}(R^\perp Y^{i}\Theta^i)\|_2^2\\
&\leq C\kappa^{-2}\|\Lambda^{\frac{1}{2}}(R^\perp Y^{i})\Theta^i\|_2^2
+C\kappa^{-2}\|R^\perp Y^{i}\Lambda^{\frac{1}{2}}\Theta^i\|_2^2\notag\\
&\leq C\kappa^{-2}\|\Lambda^{\frac{1}{2}} Y^{i}\|_2^2\|\Theta^i\|_\infty^2
+C\kappa^{-2}\|Y^{i}\|_4^2\|\Lambda^{\frac{1}{2}}\Theta^i\|_4^2\notag\\
&\leq C\kappa^{-2}\|\Lambda^{\frac{1}{2}} Y^{i}\|_2^2\|\Lambda^{\frac{3}{2}}\Theta^i\|_2^{\frac{4}{3}}\|\Theta^i\|_2^{\frac{2}{3}}
+C\kappa^{-2}\|\Lambda^{\frac{1}{2}}Y^{i}\|_2^2\|\Lambda\Theta^i\|_2^2\notag\\
&\leq \left(C\kappa^{-\frac{10}{3}}\|f\|_X^{\frac{4}{3}}M^{\frac{2}{3}}+C\kappa^{-4}\|f\|_X^2\right)
\|\Lambda^{\frac{1}{2}}Y^{i}\|_2^2\notag.
\end{align}
Applying the above estimate recursively we obtain
\begin{align}\notag
\|\Lambda^{\frac{1}{2}}Y^{i+1}\|_2^2&\leq
\left(C\kappa^{-\frac{10}{3}}\|f\|_X^{\frac{4}{3}}M^{\frac{2}{3}}+C\kappa^{-4}\|f\|_X^2\right)^i
\|\Lambda^{\frac{1}{2}}Y^1\|_2^2\\
&\leq
\left(C\kappa^{-\frac{10}{3}}\|f\|_X^{\frac{4}{3}}M^{\frac{2}{3}}+C\kappa^{-4}\|f\|_X^2\right)^i
\kappa^{-2}\|f\|_X^2\notag.
\end{align}
Therefore if we choose $\|f\|_X$ small such that 
$$
C\kappa^{-\frac{10}{3}}\|f\|_X^{\frac{4}{3}}M^{\frac{2}{3}}+C\kappa^{-4}\|f\|_X^2<1
$$
then $Y^i\to 0$  in $\dot H^{\frac{1}{2}}$ as $i\to\infty$. It implies the sequence $\left\{\Theta^i\right\}$ is Cauchy in $\dot H^{\frac{1}{2}}$ and hence has a limit $\Theta\in \dot H^{\frac{1}{2}}$. By Lemma \ref{le:l2bd} we know $\|\Theta\|_2\leq M$.

We briefly show that $\Theta$ is a solution to (\ref{SQG}) in the following. We take $i\to\infty$ in (\ref{eq:weak1}). Since $\Theta^i\to\Theta$ in $\dot H^{\frac{1}{2}}$, it follows immediately that
$$
(\Lambda^{\frac{1}{2}}\Theta^{i+1}, \Lambda^{\frac{1}{2}}\phi)\to
(\Lambda^{\frac{1}{2}}\Theta, \Lambda^{\frac{1}{2}}\phi).
$$
For the nonlinear term, we have
\begin{align}\notag
&(U^i\cdot\nabla\Theta^{i+1},\phi)-(U\cdot\nabla\Theta,\phi)\\
&=((U^i-U)\cdot\nabla\Theta^{i+1},\phi)-(U\cdot\nabla(\Theta^{i+1}-\Theta),\phi)\notag\\
&=-((U^i-U)\Theta^{i+1},\nabla\phi)+(U(\Theta^{i+1}-\Theta),\nabla\phi)\notag.
\end{align}
By the boundedness of Riesz transform and Lemma \ref{le:Sob}, it follows
\begin{align}\notag
|((U^i-U)\Theta^{i+1},\nabla\phi)|&\leq \|U^i-U\|_4\|\Theta^{i+1}\|_4\|\nabla\phi\|_2\\
&\leq \|\Lambda^{\frac{1}{2}}(\Theta^i-\Theta)\|_2\|\Lambda^{\frac{1}{2}}\Theta^{i+1}\|_2\|\nabla\phi\|_2\notag.
\end{align}
Since $\|\Lambda^{\frac{1}{2}}\Theta^i\|_2\leq C\kappa^{-1}\|f\|_X$ uniformly and $\|\nabla\phi\|_2$ is bounded, we have
$$
((U^i-U)\Theta^{i+1},\nabla\phi)\to 0 \ \ \ \mbox { as } i\to\infty.
$$
Similarly, we can show $(U(\Theta^{i+1}-\Theta),\nabla\phi)\to 0$ as $i\to\infty$.

In the end we show that $\Theta$ is the unique solution of (\ref{SQG}) among all the solutions  satisfying 
\bg\label{cond:unique}
\|\Theta\|_2\leq M, \ \ \ \|\Lambda^{\frac{1}{2}}\Theta\|_2\leq C\kappa^{-1}\|f\|_X.
\ed
Let $\tilde\Theta$ be another solution satisfying (\ref{cond:unique}). The difference $Y=\Theta-\tilde\Theta$ solves the equation
\bg\label{eq:y}
U\cdot\nabla Y+R^\perp Y\cdot\nabla\tilde\Theta+\kappa\Lambda Y=0
\ed
with $U=R^\perp Y$. Since $\Theta$ and $\tilde\Theta$ are bounded in $L^2$ and $\dot H^{\frac{1}{2}}$, we multiply (\ref{eq:y}) by $Y$ to obtain that, by applying the procedure to derive (\ref{eq:recur})
\bg\notag
\|\Lambda^{\frac{1}{2}}Y\|_2^2\leq \left(C\kappa^{-\frac{10}{3}}\|f\|_X^{\frac{4}{3}}M^{\frac{2}{3}}+C\kappa^{-4}\|f\|_X^2\right)\|\Lambda^{\frac{1}{2}}Y\|_2^2.
\ed
By the assumption on $f$ we know $\left(C\kappa^{-\frac{10}{3}}\|f\|_X^{\frac{4}{3}}M^{\frac{2}{3}}+C\kappa^{-4}\|f\|_X^2\right)<1$, which implies the solution is unique. It completes the proof of Theorem \ref{thm:ex}.

\bigskip

\section{Stability of Solutions}

In this section we investigate the stability of the steady state solutions $\Theta$ obtained in Theorem \ref{thm:ex}. Specifically, we study the Quasi-Geostrophic equation (\ref{QG}) with the initial data $\theta_0=\Theta+w_0$, where $w_0$ is considered as a perturbation. It will be established that the solutions of (\ref{QG}) approach $\Theta$ as $t\to\infty$ provided $w_0\in L^2$. Namely, we shall prove Theorem \ref{thm:stability}.

Let $w=\theta-\Theta$. By (\ref{QG}) and (\ref{SQG}) $w$ solves the equation
\begin{equation}\label{eq:equiv}
\begin{cases}
w_t+R^\perp \theta\cdot\nabla w+\kappa\Lambda w+R^\perp w\cdot\nabla\Theta=0, \\
w(0)=w_0.
\end{cases}
\end{equation}

A weak solution to \eqref{eq:equiv} is a function $w \in C_{\mathrm{w}}([0,T];L^2(\mathbb{R}^2))$ satisfying that , for any $\phi\in C_0^\infty(\mathbb R^2\times(0,T))$,
\begin{equation}\notag
\begin{split}
-\int_0^T(w,\phi_t)\, dt+\kappa\int_0^T(\Lambda^{\frac{1}{2}}w, \Lambda^{\frac{1}{2}}\phi)\, dt
&-\int_0^T(R^\perp\theta w,\nabla\phi)\, dt+\int_0^T(R^\perp w\cdot\nabla\Theta,\phi)\, dt\\
&=(w_0,\phi(x,0))+\int_0^T(f,\phi)\, dt.
\end{split}
\end{equation}

We state the existence of solutions to (\ref{SQG}) and (\ref{eq:equiv}) as follows.

\begin{Lemma}\label{thm:sqgex}
Let $\theta_0=\Theta+w_0$ with $w_0\in L^2(\mathbb R^2)$. Assume $f\in L^p$ for some $p>\frac2\alpha$. There exists a viscosity solution $\theta$ to (\ref{QG}) which satisfies, for any $T>0$,
$$
\theta\in L^\infty(0,T;L^2)\cap L^2(0,T; \dot H^{\frac{1}{2}}).
$$
Moreover, $\theta$ is bounded in $L^\infty(0,T; L^\infty(\mathbb R^2))$, that is,
\begin{equation}\notag
\th \in L^\infty(0,T; L^\infty(\mathbb R^2)).
\end{equation}
\end{Lemma}
\pf Since $\Theta\in L^2(\mathbb R^2)$ by Theorem \ref{thm:ex}, it follows $\theta_0\in L^2(\mathbb R^2)$. The existence of viscosity solutions can be obtained by stand arguments (see \cite{CC}). The boundedness of $\theta$ in $L^\infty$ on torus $\mathbb T^2$ is obtained in \cite{CD, CDmodes} (Lemma 2.3). An analogous proof  will give the boundedness on $\mathbb R^2$ (see \cite{CaV}). 
\cbdu

\begin{Lemma}\label{thm:equivex}
Let $\Theta\in H^1(\mathbb R^2)$. There exists a constant $C(\kappa)$ such that if $\|\Theta\|_{H^1}\leq C(\kappa)$, 
the equation (\ref{eq:equiv}) has a unique weak solution $w$ satisfying the energy inequality
\begin{equation}\label{energy:equiv}
\sup_{t>0}\left(\|w(t)\|_2^2+\kappa\int_0^t\|\Lambda^{\frac{1}{2}}w(s)\|_2^2ds\right)\leq \|w_0\|_2^2.
\end{equation}
\end{Lemma}
\pf
The existence can be obtained by the Galerkin approximating method, and the uniqueness follows from the linearity of the equation.
To prove the energy inequality (\ref{energy:equiv}), formally, we multiply (\ref{eq:equiv}) by $w$ and use integration by parts to infer that
\bg\notag
\frac{1}{2}\frac{d}{dt}\int_{\R^2}|w|^2dx+\kappa\int_{\R^2}|\Lambda^{\frac{1}{2}}w|^2dx
=-\int_{\R^2}R^\perp w\cdot\nabla\Theta w dx.
\ed
The right hand side can be estimated as
\begin{equation}\notag
\begin{split}
&\left|\int_{\R^2}R^\perp w\cdot\nabla\Theta w dx\right|\leq C\|R^\perp w\|_4\|w\|_4\|\nabla\Theta\|_2\\
&\leq C \|w\|_4^2\|\nabla\Theta\|_2\leq C \|\Lambda^{\frac{1}{2}}w\|_2^2\|\nabla\Theta\|_2.
\end{split}
\end{equation}
Thus if $\Theta\in H^1$ is small such that $C\|\nabla\Theta\|_2<\kappa/2$, we have
\bg\notag
\frac{d}{dt}\int_{\R^2}|w|^2dx+\kappa\int_{\R^2}|\Lambda^{\frac{1}{2}}w|^2dx\leq 0
\ed
which implies (\ref{energy:equiv}) by taking integration over the time interval $(0,t)$ and taking the limit $t\to \infty$.
\cbdu

\subsection{Proof of Theorem \ref{thm:stability}}
\label{sec-pf2}
We start with two generalized energy inequalities for the low and high frequency parts of $w$, respectively.
\begin{Lemma}\label{le:gen}
Let $\phi=e^{-|\xi|^2}, \psi=1-e^{-|\xi|^2}$, and function $E(t)\in C^1([0,\infty);L^\infty)$. Assume $\|\Theta\|_{H^1}\leq C(\kappa)$. Then the solution $w$ to (\ref{eq:equiv}) satisfies the following generalized energy inequalities,
\begin{equation}\label{gen1}
\begin{split}
\|\widecheck\phi* w(t)\|_{2}^2\leq&\|e^{-\kappa (t-s)\Lambda}\widecheck\phi* w(s)\|_{2}^2\\
&+2\int_s^t\left|(R^\perp\theta\cdot\nabla w,e^{-2\kappa (\tau-s)\Lambda}\widecheck\phi*\widecheck\phi* w(s))\right|\, d\tau\\
&+2\int_s^t\left|(R^\perp w\cdot\nabla\Theta,e^{-2\kappa (\tau-s)\Lambda}\widecheck\phi*\widecheck\phi* w(s))\right|\, d\tau,
\end{split}
\end{equation}
\begin{equation}\label{gen2}
\begin{split}
E(t)\|\psi\widehat w(t)\|_{2}^2\leq&E(s)\|\psi\widehat w(s)\|_{2}^2
-2\kappa\int_s^tE(\tau)\left\||\xi|^{\frac12}\psi\widehat w(\tau)\right\|_{2}^2\, d\tau\\
&+\int_s^tE'(\tau)\left\|\psi\widehat w(\tau)\right\|_{2}^2\, d\tau
+2\int_s^tE(\tau)\left|\left(\widehat{R^\perp w\cdot\nabla\Theta},\psi^2\widehat w\right)\right|\, d\tau\\
&-2\int_s^tE(\tau)\left|\left(\widehat{R^\perp\theta\cdot\nabla w},(1-\psi^2)\widehat w\right)\right|\, d\tau.
\end{split}
\end{equation}
\end{Lemma}
\pf We only show the a priori estimates. Multiplying (\ref{eq:equiv}) by $e^{-2\kappa (t-s)\Lambda}\widecheck\phi*\widecheck\phi* w(s)$ and integrating over the whole space 
\begin{equation}\notag
\begin{split}
\int_{\R^2}w_t e^{-2\kappa (t-s)\Lambda}\widecheck\phi*\widecheck\phi* w(s)\,dx
+\kappa\int_{\R^2}\Lambda w e^{-2\kappa (t-s)\Lambda}\widecheck\phi*\widecheck\phi*w(s)\,dx+\\
\int_{\R^2}R^\perp\theta\cdot\nabla w e^{-2\kappa (t-s)\Lambda}\widecheck\phi*\widecheck\phi* w(s)\,dx
+\int_{\R^2}R^\perp w\cdot\nabla\Theta e^{-2\kappa (t-s)\Lambda}\widecheck\phi*\widecheck\phi* w(s)\,dx=0
\end{split}
\end{equation}
The first two integrals can be written as
\begin{equation}\notag
\begin{split}
\frac12\frac{d}{dt}\int_{\R^2}\left|e^{-\kappa (t-s)\Lambda}\widecheck\phi* w(s)\right|^2\,dx
-\int_{\R^2}e^{-\kappa (t-s)\Lambda}\widecheck\phi* w(s)\partial_t\left(e^{-\kappa (t-s)\Lambda}\widecheck\phi\right)* w \,dx\\
+\kappa\int_{\R^2}e^{-\kappa (t-s)\Lambda}\widecheck\phi*w(s)\Lambda\left(e^{-\kappa (t-s)\Lambda}\widecheck\phi\right)*w \,dx\\
=\frac12\frac{d}{dt}\int_{\R^2}\left|e^{-\kappa (t-s)\Lambda}\widecheck\phi* w(s)\right|^2\,dx.
\end{split}
\end{equation}
Thus inequality (\ref{gen1}) follows from the above equation by integrating over the time interval $[s,t]$.

Taking Fourier Transform of (\ref{eq:equiv}), multiplying it by $\psi^2\widehat w E(t)$, and integrating over $\R^2$ yields
\begin{equation}\notag
\begin{split}
\frac12\frac{d}{dt}\int_{\R^2}E(t)\left|\psi\widehat w\right|^2\,d\xi-
\frac12\int_{\R^2}E'(t)\left|\psi\widehat w\right|^2\,d\xi+\kappa E(t)\int_{\R^2}|\xi|\left|\psi\widehat w\right|^2\,d\xi\\
+E(t)\int_{\R^2}\widehat{R^\perp\theta\cdot\nabla w}\psi^2\widehat w\,d\xi+E(t)\int_{\R^2}\widehat{R^\perp w\cdot\nabla\Theta}\psi^2\widehat w\,d\xi=0
\end{split}
\end{equation}
The inequality (\ref{gen2}) follows from the fact $\left(\widehat{R^\perp\theta\cdot\nabla w},\widehat w\right)=0$ and integrating the above equation on time interval $[s,t]$.
\cbdu

We use (\ref{gen1}) and (\ref{gen2}) to prove the low frequency and high frequency parts of $w$ converge to $0$ in $L^2$ respectively. The integrals on the right hand side of (\ref{gen1})-(\ref{gen2}) will be estimated in the following. 

Recall that $\theta$ is bounded in $L^\infty((0,T)\times\mathbb R^2)$, and $\|\nabla \Theta\|_2$ (hence $\|\Theta\|_\infty$) is bounded. Thus $w=\theta-\Theta$ is bounded in $L^\infty((0,T)\times\mathbb R^2)$.
Applying integration by parts, H\"older's inequality and Gagliardo-Nirenberg inequality, we have
\begin{equation}\notag
\begin{split}
&\left|(R^\perp\theta\cdot\nabla w,e^{-2\kappa (\tau-s)\Lambda}\widecheck\phi*\widecheck\phi*w(s))\right|\\
=&\left|\left(\widecheck\phi*\widecheck\phi*\Lambda^{-\frac12}\left(R^\perp\theta\cdot\nabla w\right),e^{-2\kappa (\tau-s)\Lambda}\Lambda^{\frac12}w(s)\right)\right|\\
\leq &\left\|\widecheck\phi*\widecheck\phi*\Lambda^{-\frac12}\left(R^\perp\theta\cdot\nabla w\right)\right\|_2
\left\|e^{-2\kappa (\tau-s)\Lambda}\Lambda^{\frac12}w(s)\right\|_2\\
\leq &\left\|\widecheck\phi*\widecheck\phi*\Lambda^{-\frac12}\left(R^\perp\theta\cdot\nabla w\right)\right\|_2\left\|\Lambda^{\frac12}w\right\|_2\\
\leq &\left\|\Lambda^{\frac12}\left(R^\perp\theta\right)\right\|_2\left\|w\right\|_\infty\left\|\Lambda^{\frac12}w\right\|_2
+\left\|\widecheck\phi*\widecheck\phi*R^\perp\theta\right\|_\infty\left\|\Lambda^{\frac12}w\right\|_2\left\|\Lambda^{\frac12}w\right\|_2\\
\leq &C\left\|\Lambda^{\frac12}w\right\|_2^2,
\end{split}
\end{equation}
where we used the fact that, on the low frequency part, $\widecheck\phi*\widecheck\phi*R^\perp\theta$ is bounded in $L^\infty((0,T)\times\mathbb R^2)$.
Similarly, we have
\begin{equation}\notag
\begin{split}
&\left|(R^\perp w\cdot\nabla\Theta,e^{-2\kappa (\tau-s)\Lambda}\widecheck\phi*\widecheck\phi*w(s))\right|\\
=&\left|(\widecheck\phi*\widecheck\phi*\Lambda^{-\frac12}\div(R^\perp w\Theta),e^{-2\kappa (\tau-s)\Lambda}\Lambda^{\frac12}w(s))\right|\\
\leq &\left\|\Lambda^{-\frac12}\div(R^\perp w\Theta)\right\|_2\left\|\Lambda^{\frac12}w\right\|_2\\
\leq &\left\|\Lambda^{\frac12}(R^\perp w)\Theta\right\|_2\left\|\Lambda^{\frac12}w\right\|_2+
\left\|R^\perp w\Lambda^{\frac12}\Theta\right\|_2\left\|\Lambda^{\frac12}w\right\|_2\\
\leq &\left\|\Lambda^{\frac12}(R^\perp w)\right\|_2\left\|\Theta\right\|_\infty\left\|\Lambda^{\frac12}w\right\|_2+
\left\|R^\perp w\right\|_4\left\|\Lambda^{\frac12}\Theta\right\|_4\left\|\Lambda^{\frac12}w\right\|_2\\
\leq &\left(\|\Theta\|_\infty+\left\|\nabla\Theta\right\|_2\right)\left\|\Lambda^{\frac12}w\right\|_2^2.
\end{split}
\end{equation}
Therefore, it follows that
\begin{equation}\notag
\begin{split}
\limsup_{t\to\infty}\left\|\widecheck\phi* w(t)\right\|_2^2\leq &\left(\|\Theta\|_\infty+\left\|\nabla\Theta\right\|_2+C\right)\int_s^\infty\left\|\Lambda^{\frac12}w\right\|_2^2\,d\tau\\
&+\limsup_{t\to\infty}\left\|e^{-\kappa(t-s)\Lambda}\widecheck\phi*w(s)\right\|_2^2,
\end{split}
\end{equation}
while it is known that $\limsup_{t\to\infty}\left\|e^{-\kappa(t-s)\Lambda}\widecheck\phi*w(s)\right\|_2^2=0$ since heat energy approaches $0$ as $t\to\infty$. By the energy inequality (\ref{energy:equiv}), the right hand side of the above inequality tends to $0$ as $s\to\infty$. Hence
\[\lim_{t\to \infty}\left\|\phi\widehat w(t)\right\|_2^2=\lim_{t\to \infty}\left\|\widecheck\phi*w(t)\right\|_2^2=0.\]
Next we estimate the contribution from the high frequency part. By the H\"older's inequality and Young's inequality, we have
\begin{equation}\notag
\begin{split}
\left|\left(\widehat{R^\perp w\cdot\nabla\Theta},\psi^2\widehat w\right)\right|
=&\left|\psi^2\left(\widehat{\Lambda^{-1/2}\nabla\cdot(R^\perp w\Theta)},\xi^{1/2}\widehat w\right)\right|\\
\leq &\left\|\widehat{\Lambda^{-1/2}\nabla\cdot(R^\perp w\Theta)}\right\|_2\left\||\xi|^{1/2}\widehat w\right\|_2\\
\leq &\left\|\Lambda^{\frac12}R^\perp w\right\|_2\left\|\Theta\right\|_\infty\left\|\Lambda^{\frac12} w\right\|_2+
\left\|R^\perp w\right\|_4\left\|\Lambda^{\frac12}\Theta\right\|_4\left\|\Lambda^{\frac12} w\right\|_2\\
\leq &\left(\|\Theta\|_\infty+\|\nabla\Theta\|_2\right)\left\|\Lambda^{\frac12}w\right\|_2^2
\end{split}
\end{equation}

\begin{equation}\notag
\begin{split}
\left|\left(\widehat{R^\perp \theta\cdot\nabla w},(1-\psi^2)\widehat w\right)\right|
=&\left|\left((1-\psi^2)\mathcal F(\Lambda^{-1/2}\nabla\cdot (R^\perp \theta w)), \xi^{1/2}\widehat w\right)\right|\\
\leq &\left\|(1-\psi^2)\mathcal F(\Lambda^{-1/2} \nabla\cdot(R^\perp \theta w))\right\|_2\left\||\xi|^{1/2}\widehat w\right\|_2\\
\leq &\left\|(1-\psi^2)\mathcal F(\Lambda^{1/2}(R^\perp w) w)\right\|_2\left\||\xi|^{1/2}\widehat w\right\|_2\\
&+\left\|(1-\psi^2)\mathcal F(\Lambda^{1/2}(R^\perp \Theta) w)\right\|_2\left\||\xi|^{1/2}\widehat w\right\|_2\\
&+\left\|(1-\psi^2)\mathcal F(R^\perp \theta \Lambda^{1/2} w)\right\|_2\left\||\xi|^{1/2}\widehat w\right\|_2\\
\leq &\|w\|_\infty\left\|\Lambda^{\frac12} w\right\|_2^2+\left\|\Lambda^{\frac12}R^\perp \Theta\right\|_4\|w\|_4\left\|\Lambda^{\frac12} w\right\|_2
+C\left\|\Lambda^{\frac12} w\right\|_2^2\\
\leq &(\|\nabla \Theta\|_2+C)\left\|\Lambda^{\frac12} w\right\|_2^2
\end{split}
\end{equation}
where we used the fact that, the low frequency part (associated with the factor $1-\psi^2$) of $R^\perp \theta$ is bounded in $L^\infty((0,T)\times\mathbb R^2)$, since $\theta$ is bounded in $L^\infty((0,T)\times\mathbb R^2)$.
Therefore, it follows from (\ref{gen2}) that
\begin{equation}\label{gen2new}
\begin{split}
\left\|\psi\widehat w(t)\right\|_2^2\leq &\frac{E(s)}{E(t)}\left\|\psi\widehat w(s)\right\|_2^2
-2\kappa\int_s^t\frac{E(\tau)}{E(t)}\left\||\xi|^{\frac12}\psi\widehat w(\tau)\right\|_2^2\,d\tau
+\int_s^t\frac{E'(\tau)}{E(t)}\left\|\psi\widehat w(\tau)\right\|_2^2\,d\tau\\
&+\left(\|\Theta\|_\infty+\|\nabla\Theta\|_2+C\right)\int_s^t\frac{E(\tau)}{E(t)}\left\|\Lambda^{\frac12}w(\tau)\right\|_2^2\,d\tau.
\end{split}
\end{equation}
Let $B(\rho)$ be a ball with center at origin and time dependent radius $\rho(t)$ which will be determined later. Splitting the second and third terms in (\ref{gen2new}) onto $B(\rho)$ and $B(\rho)^c$ yields
\begin{equation}\notag
\begin{split}
&-2\kappa\int_s^t\frac{E(\tau)}{E(t)}\left\||\xi|^{\frac12}\psi\widehat w(\tau)\right\|_2^2\,d\tau
+\int_s^t\frac{E'(\tau)}{E(t)}\left\|\psi\widehat w(\tau)\right\|_2^2\,d\tau\\
=&-2\kappa\int_s^t\frac{E(\tau)}{E(t)}\int_{B(\rho)}|\xi|\left|\psi\widehat w(\tau)\right|^2\,d\xi\,d\tau
-2\kappa\int_s^t\frac{E(\tau)}{E(t)}\int_{B(\rho)^c}|\xi|\left|\psi\widehat w(\tau)\right|^2\,d\xi\,d\tau\\
&+\int_s^t\frac{E'(\tau)}{E(t)}\int_{B(\rho)}\left|\psi\widehat w(\tau)\right|^2\,d\xi\,d\tau
+\int_s^t\frac{E'(\tau)}{E(t)}\int_{B(\rho)^c}\left|\psi\widehat w(\tau)\right|^2\,d\xi\,d\tau\\
\leq &\int_s^t\frac{E'(\tau)-2\kappa E(\tau)\rho(\tau)}{E(t)}\int_{B(\rho)^c}\left|\psi\widehat w(\tau)\right|^2\,d\xi\,d\tau
+\int_s^t\frac{E'(\tau)}{E(t)}\int_{B(\rho)}\left|\psi\widehat w(\tau)\right|^2\,d\xi\,d\tau.
\end{split}
\end{equation}
Choose $E(t)=(1+t)^\gamma$ and $\rho(t)=\frac{\gamma}{2\kappa (1+t)}$ with a certain large number $\gamma$ such
that
\[E'(\tau)-2\kappa E(\tau)\rho(\tau)=0.\]
Thus, it follows from (\ref{gen2new})
\begin{equation}\label{gen2end}
\begin{split}
\left\|\psi\widehat w(t)\right\|_2^2\leq &\frac{(1+s)^\gamma}{(1+t)^\gamma}\left\|\psi\widehat w(s)\right\|_2^2
+\int_s^t\frac{\gamma (1+\tau)^{\gamma-1}}{(1+t)^\gamma}\int_{B(\rho)}\left|\psi\widehat w(\tau)\right|^2\,d\xi\,d\tau\\
&+\left(\|\Theta\|_\infty+\|\nabla\Theta\|_2+C\right)\int_s^t\left\|\Lambda^{\frac12}w(\tau)\right\|_2^2\,d\tau.
\end{split}
\end{equation}
Note $\psi(\xi)\leq |\xi|^2$ for $|\xi|\leq 1$. So for large time
\begin{equation}\notag
\int_{B(\rho)}\left|\psi\widehat w(\tau)\right|^2\,d\xi\leq \int_{B(\rho)}|\xi|^4\left|\widehat w(\tau)\right|^2\,d\xi
\leq  C(1+\tau)^{-4}\| w(\tau)\|_2^2,
\end{equation}
hence, for large $s$
\begin{equation}\notag
\int_s^t\frac{\gamma (1+\tau)^\gamma}{(1+t)^\gamma}\int_{B(\rho)}\left|\psi\widehat w(\tau)\right|^2\,d\xi\,d\tau
\leq C(1+t)^{-4}\sup_{\tau\geq 0}\|w(\tau)\|_2^2
\end{equation}
which tends to $0$ as $t\to\infty$. Fix $s$ first, then
\[\frac{(1+s)^\gamma}{(1+t)^\gamma}\left\|\psi\widehat w(s)\right\|_2^2\to 0, \qquad \mbox { as } t\to\infty\]
since $\psi$ is bounded and $w\in L^\infty(L^2)$. Therefore, (\ref{gen2end}) implies
\begin{equation}\notag
\limsup_{t\to\infty}\left\|\psi\widehat w(t)\right\|_2^2
\leq \left(\|\Theta\|_\infty+\|\nabla\Theta\|_2+C\right)\int_s^\infty\left\|\Lambda^{\frac12}w(\tau)\right\|_2^2\,d\tau
\end{equation}
for arbitrary $s$. Taking $s\to\infty$ yields
\[\lim_{t\to\infty}\left\|\psi\widehat w(t)\right\|_2^2=0.\]
In the end, we conclude that
\[\lim_{t\to\infty}\|w(t)\|_2\leq\lim_{t\to\infty}\left\|\phi\widehat w(t)\right\|_2+\lim_{t\to\infty}\left\|\psi\widehat w(t)\right\|_2=0\]
which completes the proof of Theorem \ref{thm:stability}.


\medskip

\textbf{Acknowledgment.} 
The author would like to thank the anonymous referees for their valuable comments and suggestions which helped to improve the paper.


{}

\end{document}